\documentclass[11pt,a4paper]{article}

\usepackage[latin1]{inputenc}
\usepackage{amsmath}
\usepackage{amsthm}
\usepackage{amsfonts}
\usepackage{amssymb}
\usepackage{authblk}
\usepackage{indentfirst}
\usepackage{url}
\usepackage{graphicx} 
\usepackage{graphics}
\usepackage{booktabs}
\usepackage{epstopdf}
\usepackage{float}
\usepackage{caption}
\usepackage{subcaption}
\usepackage{verbatim} 
\usepackage{footnote}

\usepackage{natbib}
\usepackage[top=3cm, bottom=3cm, left=2.8cm, right=2.8cm]{geometry}
\usepackage{framed}

\usepackage{smile}

\usepackage{hyperref} % two packages for the table
\usepackage{booktabs}

\renewcommand\and{\\[\baselineskip]}

\usepackage{algorithm}
\usepackage{algorithmic}

\begin{document}
%%%%%%%%%%%%%%%%

\title{A Knowledge Gradient Policy for Sequencing Experiments to Identify the Structure of RNA Molecules Using a Sparse Additive Belief Model}

\author[1]{Yan Li\thanks{email: \texttt{yanli@princeton.edu}.}}
\author[1]{Kristofer G. Reyes\thanks{email: \texttt{kreyes@gmail.com}.}}
\author[2]{Jorge Vazquez-Anderson \thanks{email: \texttt{jorgevazquez@utexas.edu}.}}
\author[3]{Yingfei Wang\thanks{email: \texttt{yingfei@cs.princeton.edu}.}}
\author[2]{Lydia~M.~Contreras\thanks{email: \texttt{lcontrer@che.utexas.edu}. Contreras' research was supported by the Welch Foundation (Grant NO. F- 756), the Air Force Office of Scientific Research (AFOSR) Young Investigator program (FA9550-13-1-0160), and the Consejo Nacional de Ciencia y Tecnolog\'{i}a for the graduate 445 fellowship granted to J.V.A (CONACYT- 94638).}}
\author[1]{Warren B. Powell\thanks{email: \texttt{powell@princeton.edu}. Powell's esearch was supported in part by AFOSR grant contract FA9550-12-1-0200 for Program on Optimization and Discrete Mathematics, with valuable support from the program on Natural Materials and Systems.}}

\affil[1]{Department of Operations Research and Financial Engineering, Princeton University}
\affil[2]{Department of Chemical Engineering, University of Texas at Austin}
\affil[3]{Department of Computer Science, Princeton University}

\date{}
%\author{Lydia M. Contreras}
%\affil{Department of Chemical Engineering, University of Texas at Austin, Austin, Texas 78712; email: lcontrer@che.utexas.edu}
%\author
%\affil{Department of Operations Research and Financial Engineering, Princeton University, Princeton, New Jersey 08544; email: powell@princeton.edu}
\maketitle

\begin{abstract}
We present a sparse knowledge gradient (SpKG) algorithm for adaptively selecting the targeted regions within a large RNA molecule to identify which regions are most amenable to interactions with other molecules. Experimentally, such regions can be inferred from fluorescence measurements obtained by binding a complementary probe with fluorescence markers to the targeted regions. We use a biophysical model which shows that the fluorescence ratio under the log scale has a sparse linear relationship with the coefficients describing the accessibility of each nucleotide, since not all sites are accessible (due to the folding of the molecule). The SpKG algorithm uniquely combines the Bayesian ranking and selection problem with the frequentist $\ell_1$ regularized regression approach Lasso. We use this algorithm to identify the sparsity pattern of the linear model as well as sequentially decide the best regions to test before experimental budget is exhausted. Besides, we also develop two other new algorithms: batch SpKG algorithm, which generates more suggestions sequentially to run parallel experiments; and batch SpKG with a procedure which we call length mutagenesis. It dynamically adds in new alternatives, in the form of types of probes, are created by inserting, deleting or mutating nucleotides within existing probes. In simulation, we demonstrate these algorithms on the Group I intron (a mid-size RNA molecule), showing that they efficiently learn the correct sparsity pattern, identify the most accessible region, and outperform several other policies.
\end{abstract}


\section{Introduction}\label{sec::intro} 

\begin{figure*}[!htb]
\centering
\begin{tabular}{c}
\includegraphics[scale=0.6]{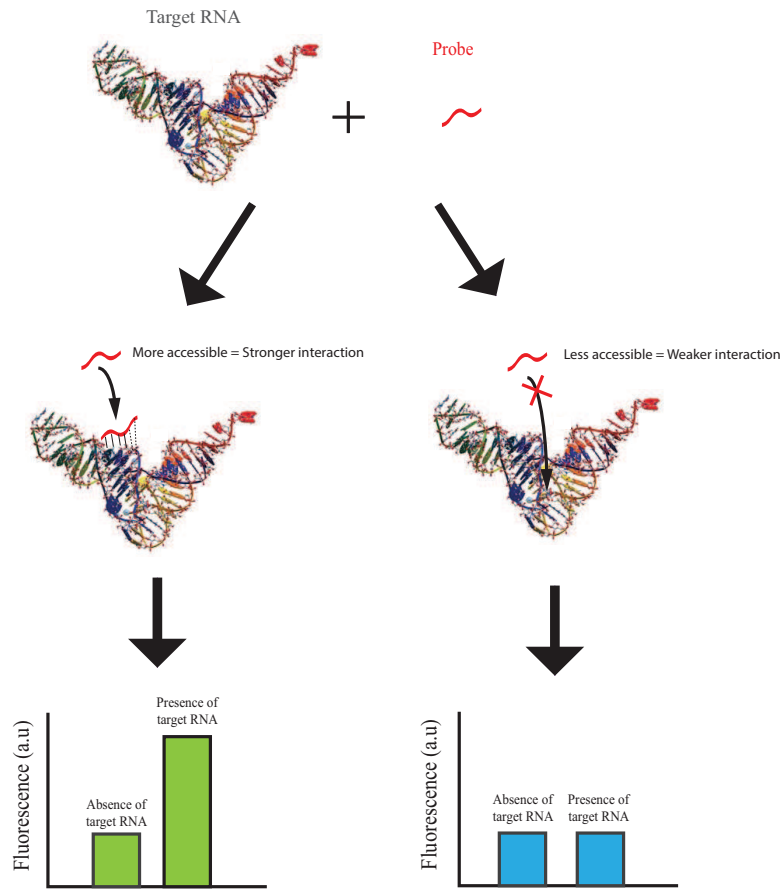}
\end{tabular}
\caption{Illustration of the in vivo RNA Structural Sensing System (iRS$^3$)}
\begin{flushleft}
 \textit{Notes.} A complementary probe with fluorescence marker is synthesized to bind to the targeted region of an RNA molecule.  If the probe targets an accessible region (left), an interaction will occur leading to a strong fluorescence signal. Otherwise, if no major fluorescence signal is detected in the presence of the target RNA, it is interpreted as no interaction between the target RNA molecule and the probe (right).
 \end{flushleft}
\label{fig:0}
\end{figure*}

In recent years RNA has been rediscovered as a potent drug target with important implications to biotechnology and human health \citep{chan2006antisense,bennett2010rna,devos2013antisense,vazquez2013regulatory,hanning2015strain}. Learning the structure of RNA molecules has become important in health research to improve the understanding of the interactions between RNA molecules and drugs. In addition, RNA regulates essential cellular processes through specific interactions with other biomolecules (e.g. proteins, other RNA or DNA molecules, etc). Disruption of an otherwise natural molecular interaction can potentially cause diseases. RNA can fold into intricate tridimensional structures making some regions accessible to interact with other molecules, while other regions remain inaccessible. Although biochemical technology has taken huge leaps by making RNA sequences readily available \citep{gelderman2013discovery}, significant progress is still required to understand RNA structure. The scientists on the team have made it possible to determine accessible regions using a fluorescence-based system \citep{sowa2015exploiting}, where an RNA strand, hereafter referred to as a probe, interacts with a specific complementary region within a target RNA generating fluorescence (see Figure \ref{fig:0}). This fluorescence directly correlates to the accessibility of a given region within a target RNA molecule. Although the in vivo RNA Structural Sensing System (iRS$^3$), as named by \cite{sowa2015exploiting}, is a valuable tool as it provides the accessibility of a segment of interest within an RNA molecule in living cells, synthesizing and running the experiments in the absence of any apriori information of the RNA can be expensive and time-consuming. With the purpose of expanding its use to characterize a full molecule, we undertake the endeavor of optimizing the experimental settings of the iRS$^3$. This paper seeks to use the knowledge gradient policy, adapted to a high-dimensional, sparse linear model to maximize the information gained from each experiment.

Our work addresses the problem of sequentially guiding experiments to identify the accessibility patterns of an RNA molecule known as the ``Tetrahymena Group I intron" (gI intron), which has been widely used as an RNA folding model, and whose complex structure has been extensively characterized \citep{cech1981vitro,kruger1982self,cech1994representation,kieft1997solution,golden1998preorganized,
russell2002exploring,koduvayur2004intracellular,wan2010multiple,vazquez2013regulatory}. Determining these accessibility patterns is difficult to do in silico, as they depend on the complicated folding of the molecule known as the tertiary structure \citep{scherr2000rna,muckstein2006thermodynamics}. For details on in silico approaches, see \cite{vazquez2013regulatory}. Experimentally, such accessibility patterns can be inferred from fluorescence measurements obtained from the iRS3 by using various complementary probes designed a priori to target a region within the intron \citep{sowa2015exploiting}. However, the number of variations of the probe increases quadratically in the number of nucleotides; therefore, the number of candidate probes is usually extremely large. A critical problem is therefore deciding which targeted regions should be tested, especially given the time and cost to perform each experiment.

The problem of identifying the accessibility pattern of an RNA molecule can be modeled mathematically as a ranking and selection (R\&S) problem. By limiting the length of the probes, we are confronted with a collection of targeted regions of the RNA, which we call alternatives in R\&S. In this R\&S problem, we have a budget of measurements that we need to allocate sequentially to test the alternatives. As more information is collected, the belief distribution is changed or updated by conditioning on all the observations we have up to this point. Our goal is to maximize our ability to gain valuable information and the reward until the budget is exhausted. Although there are many papers on R\&S problems, only a few make use of the model structure of the underlying belief model. The classical model for R\&S is the lookup table belief model, which does not assume or exploit any structure. In identifying the accessibility pattern of an RNA molecule, where we may have tens or hundreds of thousands of alternatives, the lookup table strategy is computationally intractable. It may also be inappropriate when our goal is also to learn about the underlying model structure. 

In this paper, we use a thermo-kinetic model which represents the log fluorescence level as a linear model of the weight coefficients representing the accessibility of each nucleotide. This coefficient vector is of the same dimension as that of the target molecule and thus can be high-dimensional. However, it is typically the case that only a small portion of these coefficients contain explanatory power, because not all sites are accessible due to the folding of the molecule. In such cases, a sparse linear model can offer considerably more flexibility than a linear model without sparsity structure. Therefore, we develop a sparse knowledge gradient (SpKG) policy for sequential experimental design \citep{li2015sparsekg}. This policy combines Bayesian R\&S with a frequentist learning approach for recursive Lasso (Least Absolute Shrinkage and Selection Operator) \citep{tibshirani1996regression,garrigues2008homotopy,chen2012recursive}, which is a well known $\ell_1$ regularized version of least squares. 

In this paper, we perform thorough testing of the SpKG algorithm using the setting of RNA accessibility identification.  However, the SpKG algorithm is broadly applicable to learning problems with high-dimensional belief models, a problem domain that has been attracting considerable attention.  Furthermore, SpKG can also be applied to problems with a known sparse group structure, as might happen when variables exhibit a natural clustering.  For example, applications in health might exhibit clusters of variables related to specific medical conditions.  In addition, it can also handle interaction terms with smoothing spline ANOVA models \citep[see][]{li2015sparsekg}.

It is worth noting that the SpKG algorithm is a unique and novel hybrid of Bayesian R\&S with the frequentist $\ell_1$ regularized regression known as Lasso. Both have been extensively studied in different fields \citep{friedman2001elements,powell2012optimal}. For example, as a regularized version of least squares, Lasso minimizes the residual sum of squares subject to the sum of the absolute value of the coefficients being less than a constant. As a result of the penalty term, many of the coefficients will be exactly zero for large problems. Since Lasso was first proposed by \cite{tibshirani1996regression}, there has been a considerable amount of work exploring the use of Lasso in different settings. It is useful in many settings due to its tendency to prefer solutions with fewer nonzero parameters, effectively reducing the problem dimension. For this reason, Lasso and its variants, such as elastic net regularization, are fundamental to many high-dimensional regression models. 

Most of the previous work is in the classical batch setting, where we are given a dataset from the beginning, with no control over how the observations are chosen. However, our application requires guiding experiments in an online fashion. It leads to the acquisition of new information about the environment which may improve future decisions. For this purpose, we take advantage of a homotopy algorithm for using Lasso in a recursive setting, which was proposed by \cite{garrigues2008homotopy}. This algorithm introduces an optimization problem to compute the exact update of the Lasso estimator when one additional observation is achieved. Additionally, \cite{chen2012recursive} has extended this algorithm with an $\ell_{1,\infty}$ group Lasso penalty. In this work, they consider a more general group sparsity system, which is composed of a few known nonoverlapping clusters of nonzero coefficients. The coefficients among each group have some correlation and are either all selected or not.
%In contrast to the frequentist setting, Lasso can also be interpreted as a Bayesian posterior estimate when the regression parameters have independent Laplace priors \citep{park2008bayesian}. However, since we lose conjugacy, it is infeasible to use this method for sequential learning problems. Therefore, rather than using the Bayesian Lasso, 
Our work takes Lasso estimator updated from the homotopy algorithm as a sample from the true distribution of the coefficients. Then we use this to update both the conditional normal distribution of the coefficients and the Beta-Bernoulli conjugate distribution of the probability distribution of whether each coefficient is selected or not. 

Discrete stochastic search problems have been approached in the literature from two directions: offline learning, widely known as ranking and selection (R\&S) \citep{swisher2003discrete}, and online learning, often referred to as the multiarmed bandit problem \citep{gittins2011multi}.  Either problem can be approached using frequentist or Bayesian approaches.  For example, optimal computing budget allocation, or OCBA  \citep{chen2010stochastic,chen2012optimal}, is a frequentist approach developed within the simulation optimization community for finding optimal designs that are tested using Monte Carlo simulation.  Upper confidence bounding (UCB) policies, widely studied in the machine learning community for online bandit problems, iare often approached as distribution-free strategies that enjoys bounds on the number of times that the wrong alternative might be tested \citep{auer2002finite,bubeck2012regret}.  Frequentist approaches tend to require that each alternative be tested at least once.

A substantial literature has grown up around the general strategy of maximizing the value of information. \cite{gupta1996bayesian} first proposed this idea for the offline ranking and selection problem, an idea that has been pursued under the name of the knowledge gradient (KG) \citep{frazier2008knowledge,frazier2009knowledge,powell2012optimal}.  Approximations of this idea have been proposed under names including sequential kriging optimization (SKO) \citep{huang2006global}, and efficient global optimization (EGO) \citep{jones1998efficient,bull2011convergence}.  \cite{ryzhov2012knowledge} shows that the knowledge gradient can be easily applied to both online and offline problems. Furthermore, a batch KG policy based on Monte Carlo is proposed in the decision problems where the experimentalists may be able to run several parallel experiments in batches by maximizing the value of information for an entire batch  \citep{wang2015nested}. 

Value of information policies are particularly well suited to problems where experiments are time consuming and expensive.  The idea is particularly powerful when we can exploit belief models that capture some of the underlying structure of the problem (for example, linear belief models).  This is the setting we face in this paper.

Identifying and validating RNA structures has been a problem of interest for the molecular biology community, especially since the catalytic properties of RNA were discovered back in 1982 \citep{cech1981vitro,gold1995diversity}. 
%It became apparent that those fundamental properties usually assigned to proteins by default were based upon an intimate relation between structure and function in RNA. Nowadays RNA structure studies can be counted by thousands and classified by the approaches used as experimental and computational. The former are a group of methods that includes chemical and enzymatic characterizations, most in vitro and a few in vivo \citep[see][]{sowa2015exploiting}. Although these methods provided highly valuable structural knowledge, due to their inherent complexity in general continue to be unable to catch up with the large amount of RNA sequences discovered in a daily basis. Important it is to mention that there are now a few examples of more high-throughput techniques even at a genome-wide level unveiling fundamental information on RNA structure  \citep{kertesz2010genome,ding2014vivo}. Another main group of approaches has been the use of computer-aided biophysical and stochastic predictions to the characterization of RNA structure \citep[see][for a brief review]{vazquez2013regulatory}. In this case, these approaches remain very limited to very specific types of RNAs. 
Unveiling the structure of an RNA molecule is critical to the understanding and exploitation of the interactions established with other RNA molecules. In this context structural accessibility becomes a central object of study. To this purpose the same experimental methods mentioned above have been used to understand this phenomenon, and a number of computational algorithms have been developed to identify and characterize RNA interactions \citep{pain2015assessment}. A semi-empirical thermodynamic model common to a lot of the computational algorithms available to study RNA structure is the nearest neighbor method \citep{santalucia1998unified}. A recent example of an algorithm to explain RNA-RNA interactions is the one proposed by \cite{rodrigo2013full}. This paper has proposed a thermo-kinetic model including both Gibbs free energy and a kinetic function considering an intermediate stage (known as the seeding interaction). Another important strategy to predict and understand RNA structure is the use of the partition function and other stochastic methods. In this paper, we use a novel semi-empirical (since it uses experimental DMS footprinting data) thermo-kinetic model based on the nearest neighbor model parameters \citep{xia1998thermodynamic}. We design a sequential experimentation policy based on maximizing the value of information (the knowledge gradient) using a Bayesian belief model.  We showcase this strategy in the context of the setting of characterizing the structure of an RNA molecule using fluorescent probes.

This paper makes the following contributions. (1) We provide a sparse additive belief model for the fluorescence level produced by a probe applied to an RNA molecule. Here the log fluorescence level is a linear combination of the weight coefficients describing the accessibility of each nucleotide. The coefficient vector is sparse because not all sites contribute to the thermodynamic binding process.
 (2) We derive the batch SpKG policy which generates several suggestions sequentially to run parallel experiments.
(3) We then introduce a new \textit{length mutagenesis} procedure where new alternatives, in the form of types of probes, are created by inserting, deleting, or mutating nucleotides within existing probes. Then by each experiment we enlarge the alternative library by adding in the one with the highest value of information from the larger library generated by the length mutagenesis procedure. 
 (4) We demonstrate the effectiveness of the SpKG policy, with length mutagenesis and batch learning, in the setting of selecting probes to maximize fluorescence (as an indication of identifying accessible region) for RNA molecules with hundreds of potential sites.

The remainder of the paper is organized as follows. Section \ref{sec::model} describes the ranking and selection problem and the biophysical model. In section 3, we introduce the KG policy with a (nonsparse) linear belief model, and then describe the sparse linear model proposed by \cite{li2015sparsekg}. Then we extend the SpKG algorithm to handle batch experimentation, as well as the new length mutagenesis procedure.  Section \ref{sec::simulation} reports on the application of the procedure to the in vitro DMS footprinting data with the RNA molecule Group I intron.  Section \ref{sec::conclusion} concludes the paper.

\section{Model}\label{sec::model}
We begin by considering a Bayesian R\&S model where we have $M$ alternatives. Let $\cX$ be a finite set consisting of the $M$ alternatives and $\mu_x: x \in \cX \mapsto \RR$ be a mapping from each alternative to its value. We have a budget of $N$ measurements, and we wish to sequentially decide which alternative to measure so that we can find the best alternative when our budget is exhausted. Let $\bmu = [\mu_1,\ldots,\mu_M]^T$ (Table \ref{table1} provides a summary of the notation used in the paper). We assume that $\bmu$ follows a multivariate normal distribution:
\begin{align}
\bmu \sim \cN(\btheta,\bSigma).\label{norm}
\end{align}

\begin{table}[htbp]
\caption{Table of Notation}\label{table1}
\begin{center}\footnotesize
\begin{tabular}{ll}\hline
\multicolumn{1}{l}{{\bf Variable}} & \multicolumn{1}{l}{{\bf Description}}\\ \hline
$M$ & Number of alternatives/testing probes\\
$\cX$ & Set of alternatives\\
$\mu_x$ & Unknown mean of alternative $x$\\
$N$ & Number of measurements budget\\
$\bmu$ & Column vector $(\mu_1,\ldots,\mu_M)^T$\\
$\bx^i$/$x^i$ & Sampling decision at time $i$ (vector or scalar index)\\
$y^{n+1}/y^{n+1}_x$ & Sampling observation from measuring alternative $x^n$\\
$\epsilon_{x}^{n+1}$ & Measurement error of alternative $x^n$\\
$\sigma_x$ & Known standard deviation of alternative $x$\\
$\btheta^n$, $\bSigma^n$ & Mean and covariance of prior distribution on $\bmu$ at time $n$\\
$\Pi$ & Set of all possible policies\\
$p$ & Number of features/nucleotides\\
$\phi_k(i,j)$ & Basis function for learning the local energetic value\\
$\bPhi$ & Linear transformation matrix\\
$\alpha_k$ & Weight accessibility coefficient of nucleotide at site $k$\\
$\zeta_j$ & Random indicator variable of $\alpha_j$\\
$\bvartheta^n, \bSigma^{\bvartheta,n}$ & Mean and covariance of posterior distribution on $\balpha$ after $n$ measurements\\
$(\xi_j^n, \eta_j^n)$ & Set of shape parameters of beta distribution on $p_j^n$\\
$\hat{\bvartheta}^n$ & Lasso estimate at time $n$\\
$(\hat{\bvartheta}^{n}_{\cS}, \hat{\bSigma}^{\bvartheta,n}_{\cS})$ & Mean and covariance matrix estimator from Lasso solution at time $n$\\
$S^n$ & State variable, defined as the pair $(\btheta^n,\bSigma^n)$\\
$v_x^{\text{KG},n}$ & Knowledge gradient value for alternative $x$ at time $n$\\
$L$ & Number of possible sample realizations of $\zeta^n$\\
$p_j^n$ & Parameter of Bernoulli distribution on $\zeta_j^n$\\
$K$ & Number of batch measurement budget\\
$B$ & Number of batch experiments at each time\\
$Q$ & Number of Monte Carlo simulations\\
$(\btheta^{k,b}, \bSigma^{k,b})$ & Mean and covariance of posterior distribution on $\bmu$ after $k$ batches and additional $b$ measurements\\
$(\bvartheta^{k,b}, \bSigma^{\bvartheta,k,b})$ & Mean and covariance of posterior distribution on $\balpha$ after $k$ batches and additional $b$ measurements\\
\hline
\end{tabular}
\end{center}
\end{table}

Consider a sequence of $N$ sampling decisions, $x^0,x^1,\ldots,x^{N-1}$. At time $n$, the measurement decision $x^n$ selects an alternative from set $\cX$ to sample, and we observe
\begin{align*}
y_x^{n+1} = \mu_x+\epsilon_x^{n+1},
\end{align*}
where $\epsilon_x^{n+1} \sim \cN(0, \sigma^2_{x})$, and $\sigma_x$ is known. At the beginning, we may think of $\bmu$ as a realization of the distribution given in \eqref{norm}, while the experimenter is only given some prior $\bmu \sim \cN(\btheta^0,\bSigma^0)$. Throughout the experiment, the experimenter is given the opportunity to better learn what value $\bmu$ has taken through the sequential sampling decisions. 

For convenience, let $\cF^n$ be the $\sigma$-algebra generated by the samples observed up to time $n$. Note that we have chosen our indexing so that random variables measurable with respect to the filtration $\cF^n$ are indexed by $n$ in the superscript. Following this notation, let $\btheta^n = \EE(\bmu|\cF^n)$, and $\bSigma^n = \mathrm{Var}(\bmu|\cF^n)$. This means the posterior distribution on $\bmu$ is also multivariate normal with mean $\btheta^n$ and covariance matrix $\bSigma^n$. Let $\Pi$ be the set of all $\cF^n$ measurable policies. That is $\Pi:=\{(x^0,\ldots,x^{N-1}): x^n \in \cF^n\}$. Our problem is to find the policy that solves 
\begin{align*}
\sup_{\pi \in \Pi} \EE^{\pi} \left[\max_{x \in \cX} \theta_x^N\right].
\end{align*}

\subsection{The Biophysical Model}
In the RNA accessibility identification problem, let $T$ be a molecule comprised of RNA nucleotides, called the target molecule. Denote the target molecule sequence as 
\begin{align*}
T = (t_1,\ldots,t_p),
\end{align*}
where $t_i \in \{A,C,G,U\}$ are the individual nucleotides. RNA molecules range from a few to thousands of nucleotides, but the typical length is several hundred. For the specific group I intron RNA molecule we work with in this study, it contains $p=414$ nucleotides. Depending on this sequence, the target molecule will fold upon itself in a thermodynamically favorable manner. The precise, three dimensional structure of this molecule upon folding is called the molecule's tertiary structure. Particularly, identifying regions of a molecule most amenable to interactions with other molecules is important to understanding how such interactions are mediated. Such regions depend on the molecule's tertiary structure. Those regions that are well protected in a mechanistic sense are less likely to interact with other molecules than those regions that are exposed. We refer to the regions more likely to interact with external molecules as accessible regions. Identifying such regions is accomplished by sequentially and adaptively selecting the sites of the target molecule to bind a complementary RNA probe reporter. The RNA probe reporter includes a sequence, which is typically 8 to 16 nucleotides in length. 

There is no precise definition of absolute accessibility of a region. However, we can think of the accessibility of a region relative to the accessibility of another region. This distinction is important experimentally. Therefore, to determine the accessibility of a region, the probe reporter also includes a fluorescent marker. The presence of the fluorescence at the end of an experiment, which the experimenter can measure optically, indicates whether the probe has successfully bound to the target region. The intensity of this fluorescence is an indication to how well this binding has occurred.

As described above, one can think of  an ``alternative" as a specific region within the target molecule or a complementary probe and its ``value" as the amount of binding or the fluorescence level synonymously. We now describe briefly the biophysical model that connects both through a linear model.

Our main assumption is that the fluorescence measurements are a combination of mechanistic accessibility (kinetics) and change in Gibbs free energy between bound and unbound states (thermodynamics). To model this, we consider the accessibility profile 
\begin{align*}
\balpha = (\alpha_1,\ldots,\alpha_p)^T,
\end{align*}
of the target molecule, where $\alpha_k$ is a weight describing the relative accessibility of nucleotide $k$ in the target molecule. This profile is generally unknown and can be estimated through some experimental data. Then the fluorescence intensity to target region $[i,j]$ can be modeled as 
\begin{align}
\mu_{\text{bind}}(i,j) := \log\frac{[B]}{[U]} = \phi_0(i,j) + \sum_{k=1}^{p} \alpha_k \phi_k(i,j),\label{phymodel}
\end{align}
 where $\mu_{\text{bind}} (i,j)$ represents the amount of binding to the target region $[i,j]$, $[U]$ and $[B]$ denote the fluorescence intensity of the bound and unbound states, respectively. $\phi_0(i,j)$ is a base energy gradient value for attempting to bind to region $[i,j]$, and $\phi_k(i,j)$ is the local energetic contribution of the $k$-th nucleotide position. Given the target molecule and under some assumptions, the values $\phi_k(i,j)$ are known. It is the energy as measured for Watson-Crick Helices \citep{xia1998thermodynamic}. 
 
 As mentioned above, the accessibility profile vector $\balpha = [\alpha_1,\ldots,\alpha_p]^T$ is sparse, which means that many accessible values are zero or near zero. Mechanistically, this is a reasonable assumption, as we expect the tertiary structure of any sufficiently large molecule to be well-folded, meaning the proportion of exposed, mechanistically accessible regions to protected, inaccessible regions to scale like surface area to volume. Experimentally, the prior estimate from the experimental footprinting data also shows such property as one can see later in Section \ref{sec::prior}. Therefore, this model as shown in \eqref{phymodel} is a sparse additive model. 
 
\subsection{The Bayesian Sparse Additive Model}
 The above model shows that the amount of binding $\bmu$ is linear in the weight coefficients $\balpha$ representing the accessibility of each nucleotide. Here one can view a ``feature" or an ``attribute" as the accessibility of each nucleotide. Furthermore, we let $[i^{(m)},j^{(m)}]$ represent the region for the $m$-th alternative and $\mu_m$ represent $\mu_{\text{bind}}(i^{(m)},j^{(m)})$, then we can write \eqref{phymodel} into the following affine system
\begin{eqnarray}
\begin{pmatrix}
\mu_1\\
\mu_2\\
\vdots\\
\mu_M
\end{pmatrix}
=
\begin{pmatrix}
  \phi_1(i^{(1)},j^{(1)}) & \phi_2(i^{(1)},j^{(1)})  & \cdots & \phi_p(i^{(1)},j^{(1)}) \\
  \phi_1(i^{(2)},j^{(2)}) & \phi_2(i^{(2)},j^{(2)})  & \cdots & \phi_p(i^{(2)},j^{(2)}) \\
  \vdots  & \vdots  & \ddots & \vdots  \\
  \phi_1(i^{(M)},j^{(M)}) & \phi_2(i^{(M)},j^{(M)})  & \cdots & \phi_p(i^{(M)},j^{(M)}) \\
 \end{pmatrix}
 \begin{pmatrix}
\alpha_1\\
\alpha_2\\
\vdots\\
\alpha_p
\end{pmatrix}
+
\begin{pmatrix}
\phi_0(i^{(1)},j^{(1)})\\
\phi_0(i^{(2)},j^{(2)})\\
\vdots\\
\phi_0(i^{(M)},j^{(M)})\nonumber
\end{pmatrix}.
\end{eqnarray} 
Here if we write the basis matrix as $\bPhi \in \RR^{M \times p}$ and the intercept vector as $\bPhi_0$, the above linear equations can be written in the matrix form: $\bmu = \bPhi \balpha +\bPhi_0$. Since $\bPhi_0$ is known, we can assume $\bPhi_0 = \bm{0}$ without loss of generality. Thus we have 
\begin{align}
\bmu = \bPhi \balpha, \label{lm}
\end{align} 
where $\bmu$ and $\balpha$ are random variables. We know that $\balpha$ is sparse in the sense that most of its components are zero. In the Bayesian setting, whether each feature is zero or not is also random. Specifically, let $\bzeta = [\zeta_1,\ldots,\zeta_p] \in \RR^p$ be the indicator random variable of $\balpha$, that is
 \begin{align*}
 \zeta_j = \left\{
  \begin{array}{l l}
    1 & \quad \text{if $\alpha_j \neq 0$}\\
    0 & \quad \text{if $\alpha_j=0$}
  \end{array}, \right.
  \quad \text{for } j=1,\ldots,p.
 \end{align*}
 Additionally, conditioning on $\bzeta$, we assume that $\balpha$ follows a multivariate normal distribution with mean $\bvartheta$ and covariance $\bSigma^{\bvartheta}$, that is 
 \begin{eqnarray}
\balpha \mid \bzeta \sim \cN(\bvartheta, \bSigma^{\bvartheta}). \nonumber
\end{eqnarray}
Without loss of generality, conditioning on $\bzeta$, we can permute the elements of $\balpha$ and partition $\balpha$ into the nonzero part and the zero part, so $\balpha^T = [(\balpha_{\cS})^T, \bm{0}]$. Besides, conditioning on $\bzeta$, the components of $\bvartheta$ and $\bSigma^{\bvartheta}$ are only nonzero where indexed by $\cS$. Furthermore, conditioning on $\bzeta$,  we get that $\bmu$ follows a multivariate normal distribution through the linear transformation, that is
\begin{align*}
\bmu \sim \cN(\bPhi \bvartheta, \bPhi\bSigma^{\bvartheta}\bPhi^T).
\end{align*}
Combining this with \eqref{norm}, we have
\begin{align*}
\btheta &= \bPhi \bvartheta, \\
\bSigma &= \bPhi\bSigma^{\bvartheta}\bPhi^T.
\end{align*}

The linear model in \eqref{lm} allows us to maintain the belief model in the parameter space rather than a look up table belief model in the alternative space. In the case that the parameter structure is sparse, we use a frequentist learning approach (Lasso) which uses a least squares regression with an $\ell_1$ regularization penalty, to update the belief model. In order to do this recursively, we introduce Beta-Bernoulli conjugate priors on each component of $\bzeta$.  Specifically, at time $n$, we have the following Bayesian model, for $j,j'=1,\ldots,p$,
\begin{align}
\balpha \mid \bzeta^n=\textbf{1} \sim \cN(\bvartheta^n, \bSigma^{\bvartheta,n}),\label{bayeseq1}\\
\zeta_j^n \mid p_j^n \sim \mathrm{Bernoulli}(p_j^n),\\
\zeta_j^n \perp\!\!\!\perp \zeta_{j'}^n, \quad \mathrm{for} \quad j \neq j',\\
p_j^n \mid \xi_j^n,\eta_j^n \sim \mathrm{Beta}(\xi_j^n,\eta_j^n),\label{bayeseq4}
\end{align}
where $p_j^n$ is the probability of the $j$-th feature being in the model, and $(\xi_j^n,\eta_j^n)$ are the shape parameters for the beta distribution of $p_j^n$. We assume that $\zeta_j^n$ and $\zeta_{j'}^n$ are independent for different groups $j$ and $j'$.
 Then after we get the new measurement $(\bx^n, y^{n+1})$, we recursively update the current Lasso estimate from $\hat{\bvartheta}^n$ to $\hat{\bvartheta}^{n+1}$ by the algorithm in \cite{garrigues2008homotopy} and sample a covariance matrix $\hat{\bSigma}^{\bvartheta,n+1}$ from the first order optimality condition.
 For notational simplicity, we let $\hat{\bvartheta}_{\cS}^{n+1}$ and $\hat{\bSigma}_{\cS}^{\bvartheta,n+1}$ denote the nonzero parts, leaving the $\cS$ superscripted by time $n+1$ implicit. If we regard the Lasso estimate as a sample from the conditional multivariate normal distribution, we can use the following heuristic updating scheme for a Beta-Bernoulli model and a normal-normal model. The updating equations are given by:
 \begin{align}
 \bSigma^{\bvartheta,n+1}_{\cS} = \left[(\bSigma^{\bvartheta,n}_{\cS})^{-1} + (\hat{\bSigma}^{\bvartheta,n+1}_{\cS})^{-1}\right]^{-1},\label{sigmaupdate}\\
\bvartheta^{n+1}_{\cS} = \bSigma^{\bvartheta,n+1}_{\cS}\left[ (\bSigma^{\bvartheta,n}_{\cS})^{-1} \bvartheta^{n}_{\cS}+(\hat{\bSigma}^{\bvartheta,n+1}_{\cS})^{-1} \hat{\bvartheta}^{n+1}_{\cS}\right],\\
\xi_j^{n+1}=\xi_j^n+1, \quad \eta_j^{n+1}=\eta_j^n, \quad \mathrm{for} \quad j \in \cS, \label{xiupdate1}\\
\xi_j^{n+1}=\xi_j^n, \quad \eta_j^{n+1}=\eta_j^n+1, \quad \mathrm{for} \quad j \notin \cS. \label{xiupdate2}
 \end{align}
To illustrate, one can think of the hyperparameters $(\xi_j,\eta_j)$ as the frequencies of ``in" and ``out" for each attribute and are updated through the Lasso estimates. Therefore, $\xi_j^n/(\xi_j^n+\eta_j^n)$ can be viewed as approaching the probability of the $j$-th feature being nonzero as $n$ becomes large. If the Lasso estimators can correctly recover the sparsity pattern asymptotically, then our approach should also identify the accessible nucleotides as the sampling budget goes to infinity. Theoretically, 
%previous work shows that the Irrepresentable Condition is almost necessary and sufficient for Lasso to be variable selection consistent both in the classical fixed $p$ setting and in the large $p$ setting as the sample size $N$ gets large \citep{zhao2006model}. However, several researchers have showed that such a condition does not hold in general and even if it holds \citep{zou2006adaptive}, the Lasso is not efficient for estimating the nonzero parameters. In fact, this confirms that the Lasso does not possess the oracle property \citep{fan2001variable,fan2004nonconcave}. Despite these, 
we have shown that our posterior mean estimate $\bvartheta^n$ converges to the truth $\bvartheta$ asymptotically under some conditions \citep{li2015sparsekg}.
 
\section{The SpKG Algorithms}\label{sec::KG}
Before introducing the SpKG algorithm, we first briefly review the knowledge gradient policy with both a lookup table belief model and a nonsparse, linear belief model. The knowledge gradient policy for correlated beliefs (KGCB), as introduced in \cite{frazier2009knowledge} is a fully sequential policy for learning correlated alternatives. At each time $n$, it makes the decision to measure the alternative with the largest expected incremental value, which is defined as
\begin{align*}
v_x^{\text{KG},n} = \EE (\max_{x' \in \cX} \theta_{x'}^{n+1}|S^n,x^n=x )-\max_{x' \in \cX} \theta_{x'}^n
\end{align*}
and
\begin{align*}
x^{\text{KG},n} = \arg \max_{x \in \cX} v_x^{\text{KG},n},
\end{align*}
where the knowledge state $S^n$ is defined as $S^n : = (\btheta^n,\bSigma^n)$. The KG policy can be viewed as a gradient ascent algorithm. Maintaining a multivariate normal belief on the alternative space, we can update our belief by the Bayes rule and the Sherman-Morrison formula. Taking $x^n = x$ to simplify subscripts, the updating equations are
\begin{eqnarray}
\btheta^{n+1} &=& \btheta^n+\frac{y_x^{n+1}-\theta_x^n}{\sigma_x^2+\Sigma_{xx}^n} \bSigma^n \bm{e}_x,\label{thetaupdate}\\
\bSigma^{n+1} &=& \bSigma^n - \frac{\bSigma^n \bm{e}_x \bm{e}_x^T \bSigma^n}{\sigma_x^2+\Sigma_{xx}^n}\label{lusigmaupdate},
\end{eqnarray}
where $\bm{e}_x$ is the standard basis vector with one indexed by $x$ and zeros elsewhere. We can further define a vector-valued function $\tilde{\bsigma}$ as
\begin{align*}
\tilde{\bsigma} (\bSigma^n,x) = \frac{\bSigma^n \bm{e}_x}{\sqrt{\sigma_x^2+\Sigma_{xx}^n}},
\end{align*}
and define a random variable
\begin{align*}
Z^{n+1} = \frac{(y_x^{n+1}-\theta_x^n)}{\sqrt{\mathrm{Var}[y_x^{n+1}-\theta_x^n|\cF^n]}},
\end{align*}
which is a one-dimensional standard normal random variable when conditioned on $\cF^n$ \citep{frazier2008knowledge}. Then we can write \eqref{thetaupdate} as 
\begin{align*}
\btheta^{n+1} = \btheta^n + \tilde{\bsigma} (\bSigma^n,x^n) Z^{n+1}.
\end{align*}
This allows us to compute the KG value as
\begin{align}
v_x^{\text{KG},n} &= \EE (\max_{x' \in \cX} \theta_{x'}^{n}+\tilde{\bsigma}_{x'} (\bSigma^n,x^n) Z^{n+1}|S^n,x^n=x )-\max_{x' \in \cX} \theta_{x'}^n\nonumber\\
&= h(\btheta^n, \tilde{\bsigma} (\bSigma^n,x)).\label{kgcompute}
\end{align}
Here $h: \RR^M \times \RR^M \mapsto \RR$ is defined by $h(\ba,\bb) = \EE[\max_i a_i + b_i Z]-\max_i a_i$, where $\ba$ and $\bb$ are deterministic $M$-dimensional vectors, and $Z$ is any one-dimensional standard normal random variable. \cite{frazier2009knowledge} provides a method to compute $h(\ba,\bb)$ with complexity $O(M^2 \log M)$. 

In the case of a linear model when $\bmu = \bPhi \balpha$ and we do not hold sparsity belief on $\balpha$, there exists recursive least squares (RLS) updating equations for $(\bvartheta^n, \bSigma^{\bvartheta,n})$ that are similar to the recursive updating in \eqref{thetaupdate} and \eqref{lusigmaupdate}. Before providing the updating equations, we first introduce some additional notation. Let $\bphi^n = [\phi^n_1,\ldots,\phi^n_p]^T$ be the column vector describing the alternative that is measured at time $n$ after the basis transformation $\phi$. Then, the following updating equations result from standard expressions for normal sampling of linear combinations of features \citep[see][p. 187]{powell2012optimal}:
 \begin{eqnarray}
\bvartheta^{n+1} &=& \bvartheta^n +\frac{\hat{\epsilon}^{n+1}}{\gamma^n} \bSigma^{\bvartheta,n} \bphi^n, \label{rlsmean}\\
\bSigma^{\bvartheta,n+1} &=& \bSigma^{\bvartheta,n} -\frac{1}{\gamma^n} (\bSigma^{\bvartheta,n} \bphi^n(\bphi^n)^T \bSigma^{\bvartheta,n}),\label{rlscov}
\end{eqnarray}
where $\hat{\epsilon}^{n+1}=y^{n+1}-(\bvartheta^n)^T \bphi^n$, and $\gamma^n=\sigma^2_{x}+(\bphi^n)^T \bSigma^{\bvartheta,n} \bphi^n$. For the KG computation, we can just plug the linear transformation $\btheta = \bPhi \bvartheta$ and $\bSigma= \bPhi \bSigma^{\bvartheta} \bPhi^T$ into equation \eqref{kgcompute}. The linear model exponentially reduces the computational and storage requirements of the lookup table model. The essential idea is to maintain a belief on the attributes. A summary of the KG policy with a nonsparse, linear belief model is outlined in Algorithm \ref{algorithm::KGLin}.

\begin{algorithm}
\caption{Knowledge gradient algorithm with nonsparse linear belief}\label{algorithm::KGLin}
\begin{algorithmic}[1]
%\SetKwInOut{Input}{input}\SetKwInOut{Output}{output}
%
\REQUIRE $\bvartheta^0, \bSigma^{\bvartheta,0}, \bPhi.$
 \FOR {$n=0$ to $N-1$}
   \FOR {$x = 1$ \TO $M$}
  \STATE  $\ba \leftarrow \bPhi \bvartheta^n$
  \STATE $\bb \leftarrow \bSigma^n_{x,\ast}/\sqrt{\sigma_x^2+\Sigma_{xx}^n}$
   \STATE $v \leftarrow h(\ba,\bb)$
   \IF {$x = 1$ or $v > v^{\star}$}
     \STATE $v^{\star} \leftarrow v, x^{\star} \leftarrow x$
    \ENDIF
  \ENDFOR
  \STATE{$x^n = x^{\star}$}
 \STATE Get a new measurement: $(\bx^n, y^{n+1})$ 
 \STATE $\bvartheta^{n+1} \leftarrow \bvartheta^n$ \%\textit{RLS update by \eqref{rlsmean} \eqref{rlscov}}
 \STATE $\bSigma^{\bvartheta,n+1} \leftarrow \bSigma^{\bvartheta,n}$
\ENDFOR
\RETURN $\bvartheta^N, \bSigma^{\bvartheta,N}$.
\end{algorithmic}
\end{algorithm}

\subsection{The SpKG Algorithm}\label{sec::SpKG}
It is worth noting that, for the sparse linear belief, we introduce the random variable $\zeta_j$ to indicate if the $j$-th attribute is selected or not. As a result, the KG calculation in \eqref{kgcompute} needs to be modified so that the expectation is also taken over $\bzeta$. 

Specifically, at time $n$, the Bayesian model is as described in \eqref{bayeseq1}-\eqref{bayeseq4}. The prior $\bzeta^n$ is a discrete random variable. Let $\bzeta^{n,1},\ldots,\bzeta^{n,L}$ be all the possible realizations of $\bzeta^n$, and $\PP(\bzeta^n=\bzeta^{n,l})=p^{n,l}, l=1,\ldots, L$. To compute the KG value, we need to approximate the distribution of $(\bzeta^{n+1},\bp^{n+1})$ by that of $(\bzeta^{n},\bp^{n})$. This is because the change of the sparsity belief depends on the next observation and the Lasso algorithm, and thus can be very complicated to model. Therefore, by the Law of Total Expectation, the KG value can be computed by weighting over all the possible sparsity structures \citep{li2015sparsekg}:
\begin{eqnarray}
v_x^{\text{KG},n} 
&=& \EE_{\balpha,\epsilon,\bzeta^{n+1},\bp^{n+1}} (\max_{x' \in \cX} \theta_{x'}^{n+1}|S^n,x^n=x )-\max_{x' \in \cX} \theta_{x'}^n\nonumber\\
&\approx & \EE_{\bp^n} \EE_{\bzeta^n| \bp^n} \EE_{\balpha,\epsilon|\bzeta^n,\bp^n}(\max_{x' \in \cX} \theta_{x'}^{n+1}|S^n,x^n=x, \bzeta^n,\bp^n)-\max_{x' \in \cX} \theta_{x'}^n\nonumber\\
&=& \sum^{L}_{l=1} {\EE_{p^n}(p^{n,l}) h(\ba^{n,l}, \bb^{n,l})}\nonumber\\
&=& \sum^{L}_{l=1} \prod_{\{j: \zeta_j^{n,l}=1\}} \frac{\xi_j^n}{\xi_j^n+\eta_j^n}\prod_{\{j: \zeta_j^{n,l}=0\}} \frac{\eta_j^n}{\xi_j^n+\eta_j^n} h(\ba^{n,l}, \bb^{n,l}),\label{spkgcomp}
\end{eqnarray} 
where
\begin{eqnarray}
\ba^{n,l} &=& \bPhi_{\ast,\bzeta^{n,l}} \bvartheta^{n}_{\bzeta^{n,l}},\nonumber\\
\bb^{n,l} &=& \tilde{\bsigma}(\bPhi_{\ast,\bzeta^{n,l}} \bSigma^{n,\bvartheta}_{\bzeta^{n,l}} (\bPhi_{\ast,\bzeta^{n,l}})^T, x).\nonumber
\end{eqnarray}
Here $h$ is the function defined in \eqref{kgcompute}. The subscript $\ast,\bzeta^{n,l}$ means the submatrix is taken with all the rows and columns indexed by $\bzeta^{n,l}$. The same notation is used throughout the paper. Since $L$ can be as large as $2^p$, we can sort the weights and approximate the KG value by only computing the ones with the largest probabilities. In that case, we approximately compute the KG value to avoid the curse of dimensionality. This is reasonable because we expect the sparsity pattern to converge as $n$ becomes large. We summarize the SpKG in Algorithm \ref{algorithm::KGspLin}.

\begin{algorithm}
\caption{Sparse knowledge gradient algorithm}\label{algorithm::KGspLin}
\begin{algorithmic}[1]
 \REQUIRE {$\bvartheta^0, \bSigma^{\bvartheta,0}, \{\xi_j^0, \eta_j^0\}_{j=1}^{p}, \bPhi, N$,  regularization tunable parameters $\{\lambda^i\}_{i=0}^N.$}
 \FOR {$n = 0$ to $N-1$}
 \STATE Compute KG by \eqref{spkgcomp}: $x^n = \arg \max v^{\text{KG},n}_x$ \%\textit{compute h as in Algorithm 1}
 \STATE Lasso homotopy update: $\hat{\bvartheta}^n, (\bx^n,y^{n+1}) \in \RR^m \times \RR, \lambda^n,\lambda^{n+1} \rightarrow \hat{\bvartheta}^{n+1}$
\STATE Approximately simulate $\hat{\bSigma}_{\cS}^{\bvartheta,n+1}$
\STATE $\bvartheta^{n+1} \leftarrow \bvartheta^{n}, \bSigma^{\bvartheta,n+1} \leftarrow \bSigma^{\bvartheta,n}, \{\xi_j^{n+1}, \eta_j^{n+1}\}_{j=1}^{p} \leftarrow \{\xi_j^{n}, \eta_j^{n}\}_{j=1}^{p}$ by \eqref{sigmaupdate}-\eqref{xiupdate2}
 \ENDFOR
 \RETURN {$\bvartheta^N, \bSigma^{\bvartheta,N}, \{\xi_j^N, \eta_j^N\}_{j=1}^{p}$}.
\end{algorithmic}
\end{algorithm}

\subsection{The Batch SpKG Algorithm}
The sequential knowledge gradient policy fails to account for the ability of experimentalists to run several experiments in parallel or in batches. Batch experiment means running a group of experiments at the same time. For example, the experimenter could divide a plate into squares, where each square is a different experiment, but the plate is immersed in a chemical bath at the same time. Doing experiments in parallel has the connotation of literally running different experiments on different machines at the same time. For example, in this RNA problem, the experimenter can synthesize three probes in parallel to test their fluorescence intensities all in one run. Our batch SpKG algorithm can deal with both ``batch" and ``parallel" settings. To handle such experimental settings, \cite{wang2015nested} proposes a Monte Carlo based batch knowledge gradient (BKG) approach to guide the batch experimental design by maximizing the value of information for an entire batch with the lookup table belief model. In this section, we first review the BKG policy with lookup table belief model and linear belief model, after which we then derive the new Batch SpKG policy.

To begin, we modify our notation to fit the batch measurements. Suppose we are given a batch measurement budget of $K$ with $B$ batch decisions at each time. Then the total number of measurements allowed is $N = BK$. Now at time step $k$, we choose to measure a batch of $B$ alternatives $\bx^{k,0}, \bx^{k,1}, ..., \bx^{k,B-1}$ simultaneously and get $B$ observations $y^{k+1,0},...y^{k+1,B-1}$.  We also use superscript $(k,b)$ to index our beliefs. For example, the prior multivariate normal belief can be rewritten as $(\btheta^{0,0}, \bSigma^{0,0})$. The superscript $(k,b)$ is understood as meaning that we have done $k$ batches and used $\bx^{k,0}, ... , \bx^{k, b-1}, y^{k+1,0},...y^{k+1,b-1}$ to update our belief. Similarly to equations \eqref{thetaupdate} and \eqref{lusigmaupdate}, the new updating equations can be written recursively for a batch measurements  as
\begin{eqnarray}\label{nbup}
\btheta^{k,b+1}&=& \btheta^{k,0}+ \sum_{j=0}^b\frac{y^{k+1,j}-\theta_{x^{k,j}}^{k,j}}{\sigma^2_{x^{k,j}}+\Sigma^{k,j}_{{x^{k,j}}{x^{k,j}}}}\bSigma^{k,j} e_{x^{k,j}}, \\
\bSigma^{k, b+1} &=& \bSigma^{k,b}-\frac{\bSigma^{k,b} e_{x^{k,b}}(e_{x^{k,b}})^T \bSigma^{k,b}}{\sigma^2_{x^{k,j}}+\Sigma^{k,b}_{{x^{k,b}}{x^{k,b}}}}\label{nbupp},
\end{eqnarray}
where $k=0,1,...,K-1$, $b=0,1,...,B-1$, $\btheta^{k+1,0} = \btheta^{k, B}$, and $\bSigma^{k+1,0} = \bSigma^{k, B}$. It is worth emphasizing that in the batch setting the covariance matrix would be updated within a batch since it is determined by  the measurement decisions and is independent of the observations, whereas the mean values $\btheta^{k,b}$ are only updated after the observations are collected for the whole batch. 

The batch knowledge gradient policy greedily adds in each alternative that maximizes the expected incremental value one at a time until $B$ alternatives are chosen. The expected incremental value of measuring alternatives $x_1,...,x_j$  at time step $k$ is defined as 
\begin{equation*} 
v_{x_1,... x_j}^{\text{BKG}} (S^k)= \mathbb{E}[\max_{x}\theta_{x}^{k+1}-  \max_{x} \theta^{k}_{x}| x^{k,0}=x_1,\ldots, x^{k,j-1}=x_j, S^k].
\end{equation*}
The batch knowledge gradient policy has the decision function
\begin{equation*}
x^{k,b} := X^{\text{BKG}}_b(S^{k})= \arg\max_{x \in \mathcal{X}} v_{x^{k,0},..., x^{k,b-1},x^{k,b}= x}^{\text{BKG}}(S^k),
\end{equation*}
which can be rewritten as
\begin{equation} \label{bkg}
X^{\text{BKG}}_b(S^{k})= \arg\max_{x \in \mathcal{X}}  \mathbb{E}\left[\max_{x' } \bigg(\btheta^{k, 0}+\sum_{j=0}^{b-1} \tilde{\bsigma}(\bSigma^{k,j}, x^{k,j})Z^{k+1,j} + \tilde{\bsigma}(\bSigma^{k,b},x)Z^{k+1, b} \bigg)\right],
\end{equation}
according to equation \eqref{nbup} and \eqref{nbupp}. Notice that here $x^{k,j}$, $j \le b$ are fixed when choosing $x^{k,b}$, and $\bSigma^{k,j}$ can be updated within a batch according to \eqref{nbupp}. Since an analytic expression for the expected maximization as in \eqref{bkg} is unknown, Monte Carlo sampling is used to approximate the expectation. The pseudo-code of the BKG for the $k$-th batch decision and the Monte Carlo algorithm are presented in Algorithm \ref{algorithm::Abkg} and Algorithm \ref{algorithm::mcs}. 

\begin{algorithm}
\caption{Batch knowledge gradient policy with lookup table belief for the $k$-th batch decision}\label{algorithm::Abkg}
\begin{algorithmic}[1]
\REQUIRE {$\btheta^{k,0}, \bSigma^{k,0}$, and the number of sample $Q$ for the Monte Carlo simulation}
\STATE Use the KGCB policy to find $x^{k,0}$
\STATE $\tilde{\bsigma}^0 \leftarrow \tilde{\bsigma}(\bSigma^{k,0}, x^{k,0})$
\STATE Update $\bSigma^{b,1}$ according to \eqref{nbupp}
\FOR{$b=1$ to $B-1$}
 \STATE Use Algorithm \ref{algorithm::mcs} below to calculate $v_{x^{k,0},... x^{k,b-1},x^{k,b}=x}^{\text{BKG}}$
 \STATE $x^{k,b}= \arg\max_{x \in \mathcal{X}} v_{x^{k,0},..., x^{k,b-1},x^{k,b}= x}^{\text{BKG}}$
 \STATE $\tilde{\bsigma}^b \leftarrow \tilde{\bsigma}(\bSigma^{n,b}, x^{n,b})$
 \STATE Update $\bSigma^{n,b+1}$ according to \eqref{nbupp}
\ENDFOR
\RETURN {batch decisions $x^{k,0}, x^{k,1},..., x^{k,B-1}$}
\end{algorithmic}
\end{algorithm}

\begin{algorithm}
\caption{Monte Carlo simulation for calculating KG values}\label{algorithm::mcs}
\begin{algorithmic}[1]
\REQUIRE{$b$, $\btheta^{k,0}, \tilde{\bsigma}^0, \tilde{\bsigma}^1, ..., \tilde{\bsigma}^{b-1}, \bSigma^{k, b}$, and $Q$}
\FORALL{$x \in \mathcal{X}$}
\STATE $\text{sum}_{x}=0$
\FOR{$q=1$ to $Q$}
   \FOR{$j=0$ to $b$}
     \STATE     Generate a realization $z^j_q$ of $Z^{k,j}$
   \ENDFOR
  \STATE    temp  $\leftarrow \max_{x'}\big (\btheta^{k,0}_{x'}+ \sum_{j=0}^{b-1}\tilde{\bsigma}^j_{x'}z^j_q+\tilde{\bsigma}(\bSigma^{k,b}, x)z^b_q \big)$
\STATE        sum$_{x}$ $\leftarrow$ sum$_{x}+$ temp
\ENDFOR
\ENDFOR
\RETURN{\textbf{sum}/Q }
\end{algorithmic}
\end{algorithm}

Besides, the logic of the batch knowledge gradient policy can be generalized to a linear belief model. In this case, instead of recursively updating $\btheta^{k,b}$ and $\bSigma^{k,b}$ directly as in \eqref{nbup} and \eqref{nbupp}, we recursively update  $\bvartheta^{k,b}$ and $\bSigma^{\bvartheta,k,b}$ for a batch of observations through RLS and use the linear transformation $\btheta = \bPhi \bvartheta$ and $\bSigma= \bPhi \bSigma^{\bvartheta} \bPhi^T$:
 \begin{eqnarray}
\bvartheta^{k+1,b+1} &=& \bvartheta^{k,0} +\sum_{j=0}^b\frac{\hat{\epsilon}^{k+1,j}}{\gamma^{k,j}} \bSigma^{\bvartheta,k,j} \bphi^{k,j},\label{btheta}\\
\bSigma^{\bvartheta,k+1,b+1} &=& \bSigma^{\bvartheta,k,b} -\frac{1}{\gamma^{k,b}} (\bSigma^{\bvartheta,k,b} \bphi^{k,b}(\bphi^{k,b})^T \bSigma^{\bvartheta,{k,b}}), \label{bsigma}
\end{eqnarray}
where $\hat{\epsilon}^{k+1,j}=y^{k+1,j}-(\bvartheta^{k,j})^T \bphi^{k,j}$, and $\gamma^{k,j}=\sigma^2_{x}+(\bphi^{k,j})^T \bSigma^{\bvartheta,k,j} \bphi^{k,j}$.

Furthermore, for the sparse linear model, the posterior distribution of the sparsity structure parameter $\bzeta$ will be updated according to \eqref{xiupdate1} and \eqref{xiupdate2} after the observations are revealed for the whole batch.  Since $\balpha|\bzeta^k \sim \cN(\bvartheta^k, \bSigma^{\bvartheta,k})$, then given $\bzeta$, $\bvartheta^{k,b}_{\bzeta^{k,l}}$ and $\bSigma^{\bvartheta,k,b}_{\bzeta^{k,l}}$ can be updated according to \eqref{btheta} and \eqref{bsigma}. The batch SpKG algorithm works by greedily adding in each alternative that maximizes the expected marginal value until $B$ alternatives are chosen given the sparsity structure unchanged within a batch. The KG value of measuring alternatives $x_1,...,x_j$ at time step $k$ can be computed as 
\begin{equation}\label{BS}
v_{x^{k,0},... x^{k,b-1},x^{k,b}=x}^{\text{BSpKG}}(S^k)
= \sum^L_{l=1} {\EE_{p^k}(p^{k,l})}  v_{x^{k,0},... x^{k,b-1},x^{k,b}=x}^{\text{BKG}}( S^{\text{Sp},k,l}),
\end{equation}
where $S^{\text{Sp},k,l}$ is the knowledge state given $\bzeta^{k,l}$, and $S^{\text{Sp},k,l}=\big(\bPhi_{\ast,\bzeta^{k,l}} \bvartheta^{k,0}_{\bzeta^{k,l}}, \bPhi_{\ast,\bzeta^{k,l}} \bSigma^{\bvartheta,k,0}_{\bzeta^{k,l}} (\bPhi_{\ast,\bzeta^{k,l}})^T\big)$. After $x^{k,0},...,x^{k,b-1}$ are chosen,  $v_{x^{k,0},... x^{k,b-1},x^{k,b}=x}^{\text{BKG}}( S^{\text{Sp},k,l})$ can be approximated using Algorithm \ref{algorithm::mcs}. Then the $b$-th decision within a batch is given by $x^{k,b}=\arg\max_{x}v_{x^{k,0},... x^{k,b-1},x^{k,b}=x}^{\text{BSpKG}}(S^k)$.

The batch SpKG algorithm can be summarized in Algorithm \ref{algorithm::BSpKG}.
\begin{algorithm}
\caption{Batch sparse knowledge gradient policy for the $k$-th batch decision}\label{algorithm::BSpKG}
\begin{algorithmic}[1]
 \REQUIRE{$\bvartheta^{k,0},  \bSigma^{\bvartheta,k,0}$, and the number of sample $Q$ for the Monte Carlo simulation}
 \STATE Use Algorithm \ref{algorithm::KGspLin} to find $x^{k,0}$;\\
   \FOR{$l=1$ to $L$}
      \STATE $\tilde{\bsigma}^{0,l} \leftarrow \tilde{\bsigma}\big(\bPhi_{\ast,\bzeta^{k,l}} \bSigma^{\bvartheta,k,0}_{\bzeta^{k,l}} 		(\bPhi_{\ast,\bzeta^{k,l}})^T, x^{k,0}\big)$
      \STATE Update $\bSigma^{\bvartheta,k,1}_{\bzeta^{k,l}}$ according to \eqref{bsigma}
    \ENDFOR
\FOR {$b=1$ to $B-1$}
	\FOR{$l=1$ to $L$}
      \STATE	Use Algorithm \ref{algorithm::mcs} to calculate $v_{x^{k,0},... x^{k,b-1},x^{k,b}=x}^{\text{BKG}}(S^{\text{Sp},k,l})$ with input parameters $\btheta^{k,0} =\bPhi_{\ast,\bzeta^{k,l}} \bvartheta^{k,0}_{\bzeta^{k,l}}$, $\bSigma^{k,b}= \bPhi_{\ast,\bzeta^{k,l}} \bSigma^{\bvartheta,k,b}_{\bzeta^{k,l}} (\bPhi_{\ast,\bzeta^{k,l}})^T$, and $\tilde{\bsigma}^{0,l},...,\tilde{\bsigma}^{b-1,l}$
    \ENDFOR
  \STATE Calculate $v_{x^{k,0},... x^{k,b-1},x^{k,b}=x}^{\text{BSpKG}}$ according to \eqref{BS}
  \STATE $x^{k,b}= \arg\max_{x \in \mathcal{X}} v_{x^{k,0},..., x^{k,b-1},x^{k,b}= x}^{\text{BSpKG}}$
\FOR{$l=1$ to $L$}
  \STATE $\tilde{\bsigma}^{b,l} \leftarrow \tilde{\bsigma}\big(\bPhi_{\ast,\bzeta^{k,l}} \bSigma^{\bvartheta,k,b}_{\bzeta^{k,l}} 		(\bPhi_{\ast,\bzeta^{k,l}})^T, x^{k,b}\big)$
  \STATE Update $\bSigma^{\bvartheta,k,b+1}_{\bzeta^{k,l}}$ according to \eqref{bsigma}
\ENDFOR
\ENDFOR
\RETURN {batch decisions $x^{k,0}, x^{k,1},..., x^{k,B-1}$}
\end{algorithmic}
\end{algorithm}

\subsection{The Batch SpKG Algorithm with Length Mutagenesis}
For all of the above algorithms, we use a fixed set of discrete alternatives $\cX$. For this RNA accessibility identification problem, if we consider the probe sequence with length of 8$\sim$16, the number of alternatives can be $\sim$4000. Including all of them in the alternative library would be computationally expensive while working with only a small subset would possibly miss the most accessible region. To compromise, we propose a novel procedure that we call \textit{length mutagenesis}, which sequentially enlarges the probe library through an adaptive probe refinement procedure. We simply refer to this as ``mutagenesis" throughout the remainder of the paper for compactness. The mutagenesis works as follows. Suppose that at time $n$ we have a library of probes that we denote as
\begin{align*}
\cX^n = \{x_1,\ldots,x_{M^n}\}.
\end{align*}
For the next experiment, we could consider a larger library of sub-probes obtained through mutagenesis. With mutagenesis, we either add or delete nucleotides at one end of a probe. For now, let us think of each alternative $x$ as a probe sequence representing a region in the target molecule. Specifically, given a probe $x$, we can alter it through a round of mutagenesis to get a new probe $x'$ of the form
\[
x = [i, j] \rightarrow x' = \left\{
  \begin{array}{l l}
   { [i+k, j]}, & \quad  i+k < j\\
   { [i+k, j]}, & \quad i < j+k
  \end{array}, \right.
 0< |k| \leq 7.
\]
Since a probe of length less than 4 does not necessarily bind to the correct targeted region experimentally, we limit the probe length to be no less than 4. Then for a probe $x$, let $\cM(x)$ denote the set of possible probes obtained from $x$ through mutagenesis. At time $n$, we get an expanded library through mutagenesis, that is
\begin{align*}
\bar{\cX}^n = \cX^n \cup \bigcup_{i = 1}^{M^n} \cM(x_i).
\end{align*}
From this expanded library, we pick the alternative with the highest KG score, that is
\begin{align*}
x^n = \arg\max_{x \in \bar{\cX}^n} v_x^{\text{KG},n}.
\end{align*}
Then we add this probe to our library if it is new, 
\begin{align*}
\cX^{n+1} = \cX^n \cup \{ x^n\}.
\end{align*}
This approach allows us to add a new probe which potentially has a higher fluoresence signal at each time and work dynamically with the alternative library. The computation is simpler and more efficient than maintaining all the possible alternatives. As shown in the simulations in Figure \ref{fig:8}, we have a much better chance of obtaining information about the accessibility of the molecule relative to maintaining a fixed probe region. %The batch SpKG algorithm with this adaptive probe refinement procedure is summarized in Algorithm .

\section{Empirical Study}\label{sec::simulation}
In this section, we present the simulation results of the SpKG algorithms for the RNA accessibility identification problem described in Section \ref{sec::model}. The target RNA molecule is the Group I intron, a mid-size RNA molecule sequence with a length of 414. Due to the nature of this problem, we are not able to identify the 21 nucleotides at one end of the molecule sequence. Therefore, we only work with a sequence with length of 393. 

In Section \ref{sec::prior}, we describe the prior in vitro DMS footprinting data and the methods for generating the prior covariance matrix. In Section \ref{sec::simres}, we present the simulation results of the performance of SpKG on a collection of probe sequences  as well as those of the batch SpKG algorithm and the batch SpKG with mutagenesis scheme. We also compare these policies with several other policies. 

\subsection{Prior Distribution}\label{sec::prior}
When choosing a prior distribution, the domain experts can have many ways to articulate their prior belief on the accessible regions.  In this problem, we have the accessibility profile obtained from the in vitro DMS footprinting. The DMS footprinting is a standard chemical method to study RNA structure. It relies on the reactivity of a small molecule Di-methyl sulfate (DMS) with the base-pairing molecular faces of adenosines and cytidines (A and C). The higher the DMS reactivity is for a nucleotide site, the more the nucleotide is exposed. By reversely transcribing the DMS reacted RNA into DNA, we can determine sites of reaction and thus the levels of protection exposure at a single-nucleotide resolution. Here we use in vitro DMS data from \cite{russell2006paradoxical} as an initial estimation of nucleotide accessibility. One may think of this dataset as providing the priors $\bvartheta^0$ and $(\xi_j^0, \eta_j^0)$ for $j=1,\ldots,p$. 

\begin{figure*}[!htb]
\centering
\begin{tabular}{c}
\includegraphics[scale=0.5]{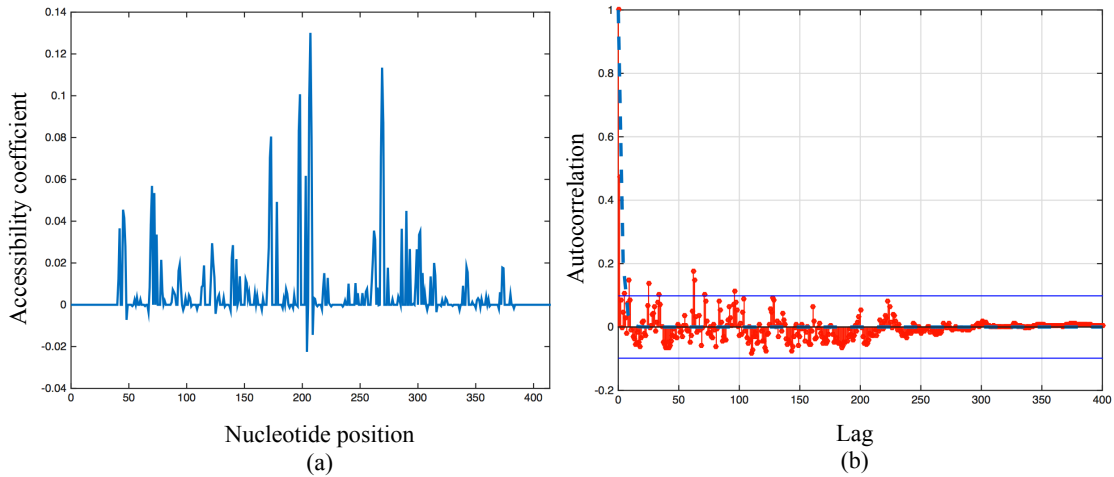}
\end{tabular}
\caption{(a) Prior Accessibility Coefficients of the In vitro DMS footprinting Data; (b) The Sample Autocorrelation Function for Prior Accessibility Coefficients}
\begin{flushleft}
 \textit{Note.} The fitted exponential decay function with a decay rate of 0.39728 is plotted in the blue dash line in (b).
 \end{flushleft}
\label{fig:1}
\end{figure*}

We now discuss how we generate the prior covariance matrix $\bSigma^{\bvartheta,0}$. For some previous work, this matrix is generated by taking the diagonal matrix with the variance from the measurement noise. This means we begin with independent beliefs. However, for this problem, the weight accessibility coefficients have natural proximity correlations. As can be seen from Figure \ref{fig:1} (a), the value of the accessibility coefficients are quite close locally. In fact, if we plot the sample autocorrelation function as shown in Figure \ref{fig:1}(b), we can see that the correlation is 0.4718 when the lag is 1, jumps in the interval [-0.1, 0.2] for lag until 100 and almost decays to 0 after 250. If we fit an exponential function $y=e^{-\kappa x}$ to the sample autocorrelation, the best fitted decay rate by least squares is $\kappa^{\star} = 0.39728$. 

This proximity correlation detection is also consistent with the domain experts' experience that the probes tend to perform similarly within a window of up to $\sim$40. Theoretically, the larger this window is, the less measurements we need to identify the accessibility patterns of the target RNA. This is because the advantage of our algorithm is to incorporate the covariance matrix $\bSigma^{\bvartheta}$, so we can locally infer more information based on what we have learned. Taking advantage of this proximity correlation knowledge, we use the exponential covariance function to model the prior covariance $\bSigma^{\bvartheta,0}$. Under the exponential covariance functions, for any two points $i$ and $j$ from $1,\ldots,p$,
\begin{align}
\mathrm{Corr}(\vartheta_i, \vartheta_j) &=  \exp\{ -\gamma|i-j |\},\label{corr}\\
\mathrm{Var}(\vartheta_i) &= \beta_i^2, \label{var}
\end{align}
where $\gamma >0$ and $\beta_1,\ldots,\beta_p >0$ are hyperparameters chosen to reflect our belief. Here, $\beta_i$ should be chosen to represent our confidence that $\vartheta_i$ is close to our chosen mean function. The $\gamma$ should reflect how quickly we believe $\Sigma^{\bvartheta}_{i,j}$ changes as $i$ and $j$ move further apart, with larger values of $\gamma$ suggesting more rapid change. This simple family of covariance functions produces Gaussian process priors that are stationary and thus can be used for modeling the accessibility coefficients in this problem. 

In practice, when one is unsure about the value of these hyperparameters, second-level priors can be put to model the coefficients $\gamma$ and $\beta_i$. However, instead of using these hierarchical maximum a posteriori (MAP) approaches, we directly set up the values according to the prior in vitro DMS footprinting data and the domain experts' experience. Specifically, we let $\gamma = \kappa^{\star}=0.39728$ from the fitted decay rate of the in vitro DMS footprinting data as shown in Figure \ref{fig:1}(b). Also, according to our domain experts, the noise ratio for estimating the accessibility coefficient is $10\% \sim 15\%$, so we set the noise ratio  $r = 20\%$ to be conservative. That is
\begin{align}\label{noise}
\beta_j = 20\% \times \tilde{\vartheta}_j, \text{ for } j = 1,\ldots,p,
\end{align}
where 
\begin{align}\label{tilvarth}
\tilde{\vartheta}_j := \left\{
\begin{array}{l l}
 \vartheta_j & \text{ for } j: \vartheta_j \neq 0\\
 \sum_{j=1}^p \vartheta_j/p & \text{ for } j: \vartheta_j = 0,
\end{array} \right.
\end{align}
Combining \eqref{corr}, \eqref{var}, and \eqref{noise}, we set the prior covariance matrix $\bSigma^{\bvartheta,0}$ as
\begin{align}\label{pcov}
\Sigma^{\bvartheta,0}_{i,j} = r^2 \tilde{\vartheta}_i \tilde{\vartheta}_j \exp\{-\kappa^{\star} |i-j | \}, \quad \text{for } i,j = 1,\ldots,p,
\end{align}
where $r = 20\%$, $\kappa^{\star} = 0.39728$, and $\tilde{\vartheta}_1,\ldots\tilde{\vartheta}_p$ are from \eqref{tilvarth}.

Besides the prior covariance matrix, we also have to set the shape parameters $(\xi_j^0,\eta_j^0)$ for the beta distribution (the frequency priors) in \eqref{bayeseq4}. For $j=1\ldots,p$, we propose to set the frequency priors as 
\begin{align*}
(\xi_j, \eta_j) =\left\{
\begin{array}{l l}
(1, 1) + (w,0), & \text{ for } j: \vartheta_j \neq 0\\
(1, 1) + (0, w), & \text{ for } j: \vartheta_j = 0
\end{array}. \right.
\end{align*}
Here $w \geq 0$ is a hyperparameter representing our confidence in the prior sparsity pattern. A smaller $w$ reflects less confidence in the prior while a larger $w$ represents more confidence. If at the end of the experiments our algorithm uses probability 0.5 as a threshold to choose the nonzero coefficients, then $w$ should not be larger than the sampling budget $N$. Otherwise, the sparsity pattern of the posterior estimate is totally identical with the prior data no matter what Lasso estimates we get. In the following simulations, we treat $w$ as a tunable parameter depending on either the good prior or the bad prior cases.

\subsection{Simulation Results}\label{sec::simres}
Notice that for this RNA accessibility identification problem, the real accessibility profile is unknown, while we can approximately learn this through various experimental methods. In this paper, we perform various simulations in which we sample a truth from a stochastic process and then run this trial with this fixed truth for some fixed number of measurement budget $N$. Then we replicate this over several runs to assess the performance of various policies. This truth coefficient is usually sampled through both vertically perturbing the values of the prior coefficient by a normally distributed random variable, and horizontally rotating the prior along the RNA molecule. In this section, we show the simulation results for all of the algorithms demonstrated in Section \ref{sec::KG}.

\subsubsection{Results using the SpKG Algorithm}
We begin by describing the results for the SpKG policy presented in Algorithm \ref{algorithm::KGspLin}. In this simulation, as suggested by the scientists on the team, we try all the probes of length $10$ with $3$ overlaps for the adjacent ones. Then we have $M = 55$ number of alternatives. In this setting, we compare Algorithm \ref{algorithm::KGspLin} with two other policies: a pure exploration policy, which randomly chooses an alternative to test at each time, and Algorithm \ref{algorithm::KGLin}, which uses KG with a nonsparse, linear belief model. It is worth emphasizing that for the pure exploration policy here, although it does not use the sparse linear structure to make measurement decisions, it still updates the belief in the same way as that of SpKG. 

In this simulation, we generate the true accessibility coefficient vector from a multivariate normal distribution with mean $\bvartheta^0$, with the covariance matrix being the same form as \eqref{pcov}. The differences are the noise ratio $r$ is chosen to be as large as $10$, and $\kappa$ is drawn from a normal distribution with mean $\kappa^{\star}$ and standard deviation 0.1. Then we take this vector and right circularly shift it by an integer value uniformly sampled from 20 to 50 at each time. Since the prior has now been significantly altered from the truth, we believe it is not a good prior and set $w = 10$ with the measurement budget $N=100$. To quantitatively measure the performance of different policies, we consider the opportunity cost (OC), defined as the difference between the true value of the alternative that is actually the best and the true value of the alternative that is the best according to the policy's posterior belief distribution, i.e.,
\begin{align*}
\text{OC}(n) = \mu(x^{\star})-\mu(x^{n,\star}),
\end{align*}
where $x^{\star}$ is the true optimal alternative, and $x^{n,\star}$ is the estimated optimal alternative at time $n$. So OC describes how far from optimal the current estimate of the optimal solution is after each experiment and and thus can serve as a metric for the performance of a specific decision policy. For illustrative purposes, we also consider the percentage OC with respect to the optimal value,
\begin{align*}
\text{OC}\%(n) = \frac{\mu(x^{\star})-\mu(x^{n,\star})}{\mu(x^{\star})}.
\end{align*}
This normalized representation is unit-free and better illustrates how far in percentage we are from the optimal. 

\begin{figure*}[!htb]
\centering
\begin{tabular}{c}
\includegraphics[scale=0.55]{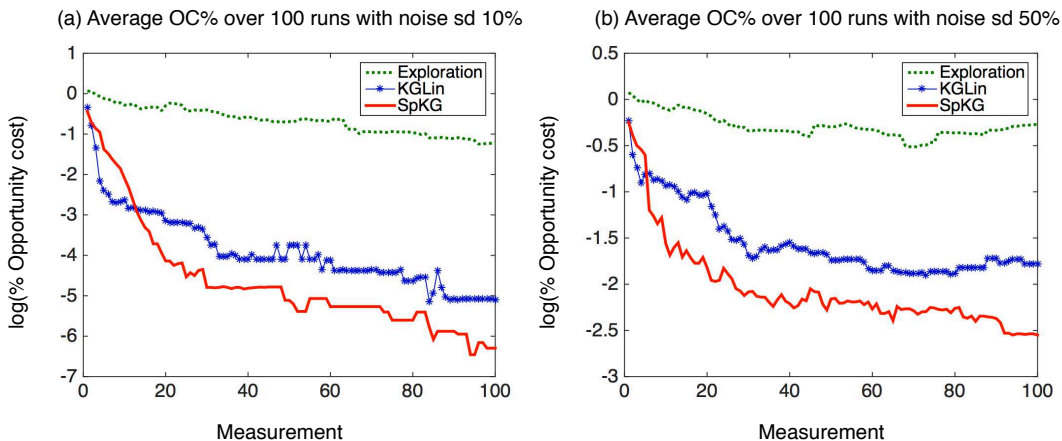}
\end{tabular}
\caption{Average Opportunity Cost Comparing Policies with Low and High Measurement Noise}
\begin{flushleft}
\textit{Notes.} This simulation is for the whole target molecule sequence. The alternative probes are of length $10$ with $3$ overlaps for the adjacent ones. To better visualize the difference in OC, the average percentage OC is plotted on a log scale over 100 simulation trials.
\end{flushleft}
\label{fig:2}
\end{figure*}

Figure \ref{fig:2} plots the average percentage OC on a log scale over 100 runs for three different policies in both low and high noise settings. The standard deviations of noise are 10\% and 50\% of the expected range of the truth. As both figures show, compared with pure exploration, SpKG results in a significant reduction in the opportunity cost, although the exploration policy also takes advantage of the sparse linear structure for the Bayesian implementation. When comparing with KG for a nonsparse, linear belief model, SpKG outperforms with lower average opportunity costs most of the time. However, during the initial stage when there are less than five measurements, SpKG behaves no better than KGLin, especially in low noise settings. This is because it takes several samples for Lasso to identify the true support. When it finds the sparsity pattern, SpKG is able to find alternatives converging to the truth much more efficiently than other policies.

\subsubsection{Results using the Batch SpKG Algorithm}\label{sec::simBKG}
In this simulation, we try testing the batch SpKG algorithm described in Algorithm \ref{algorithm::BSpKG}. For the real experiments, the experimentalist is able to synthesize three probes to test the fluorescence intensities in parallel at each time. So in the batch setting, we let $B = 3$. From now on, let us take a specific region of the RNA molecule from site 95 to site 251, as this region is suggested by the domain experts to be the most promising sub-sequence. We try a larger set of probes than before: all 8-long probes shifted by 4, all 12-long probes shifted by 6, and all 16-long probes shifted by 8. That is we take all the regions $[4k+1, 4k+8], [6k+1,6k+12]$, and $[8k+1, 8k+16]$ starting from 95 and ending at 251. Beside these, we also include ten other probes as suggested by the domain experts: $[98, 112]$, $ [113, 126]$, $ [127, 140]$, $ [141, 155]$, $[156, 170]$, $[171,179]$, $[179,194]$, $[195,214]$, $[215,233]$, and $[234,251]$. In total, we have $M = 91$ number of alternatives. We generate the true accessibility profile in the same way as before. 

\begin{figure*}[!htb]
\centering
\begin{tabular}{c}
\includegraphics[width=0.95\linewidth]{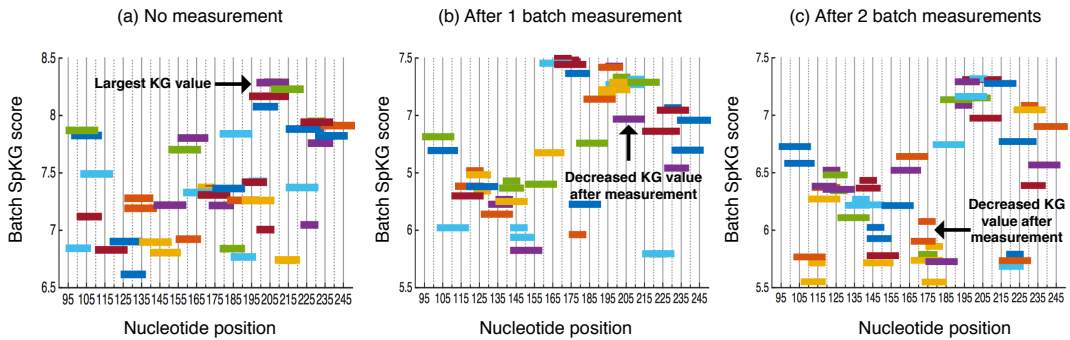}
\end{tabular}
\caption{Batch SpKG Values Before and After 1 and 2 Batch Measurements with Noise Ratio of 30\%. }
\begin{flushleft}
\textit{Notes.} This simulation is for a selected set of probes ranging from site 95 to site 251. Each bar is a potential range of a probe. The vertical axis is the batch SpKG score. Only those probes with the batch SpKG score above the mean are plotted. The arrows indicate the decreases in KG values for the probes that were previously measured. Note that some of these are not shown since they have KG values below average.
\end{flushleft}
\label{fig:5}
\end{figure*}

First, we illustrate how the batch SpKG policy works under a measurement noise of $30\%$. For one such simulated truth, we depict the batch SpKG value initially, after one, and two batch measurements, respectively in Figure \ref{fig:5}. For these figures, we only include those probes with batch SpKG values above the mean to better visualize the KG scores. As indicated by the arrows, for the probes with the largest batch SpKG scores, the KG scores drop after they have been measured. As we only plot those with KG scores above average, some probes with high KG scores in Figure \ref{fig:5}(a) have the scores dropped below average after being measured and are therefore not shown in Figure \ref{fig:5}(b). %Specifically, as suggested by the Batch SpKG policy,  the best triplets to measure are probes $[199, 214]$, $[195, 202]$ and $[167, 178]$. The first one is indicated in Figure \ref{fig:5}(b) with a dropped KG score 
This observation is also consistent with our intuition of SpKG as a measure of the value of information, and thus we can use this policy as a guideline to pick the next experiments.

\begin{figure*}[!htb]
\centering
\begin{tabular}{c}
\includegraphics[scale=0.5]{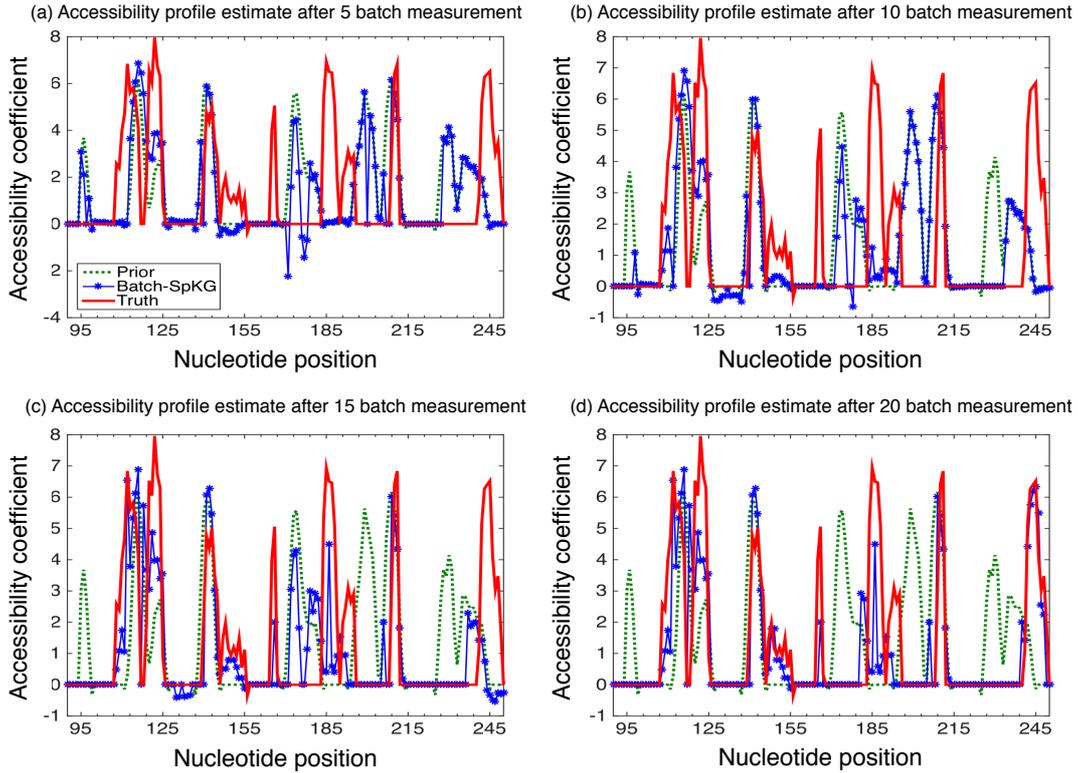}
\end{tabular}
\caption{Accessibility Profile Estimate by the Batch SpKG Algorithm After 5, 10, 15, and 20 Batch Measurements with Noise Ratio of 20\%}
\begin{flushleft}
\textit{Note.} This is for the accessibility coefficient ranging from site 95 to site 251. 
\end{flushleft}
\label{fig:3}
\end{figure*}

Furthermore, for one simulated truth, we also plot the estimates of accessibility profiles (coefficients) after 5, 10, 15, and 20 batch measurements with a noise ratio of 20\% in Figure \ref{fig:3}.

\begin{figure*}[!htb]
\centering
\begin{tabular}{c}
\includegraphics[scale=0.65]{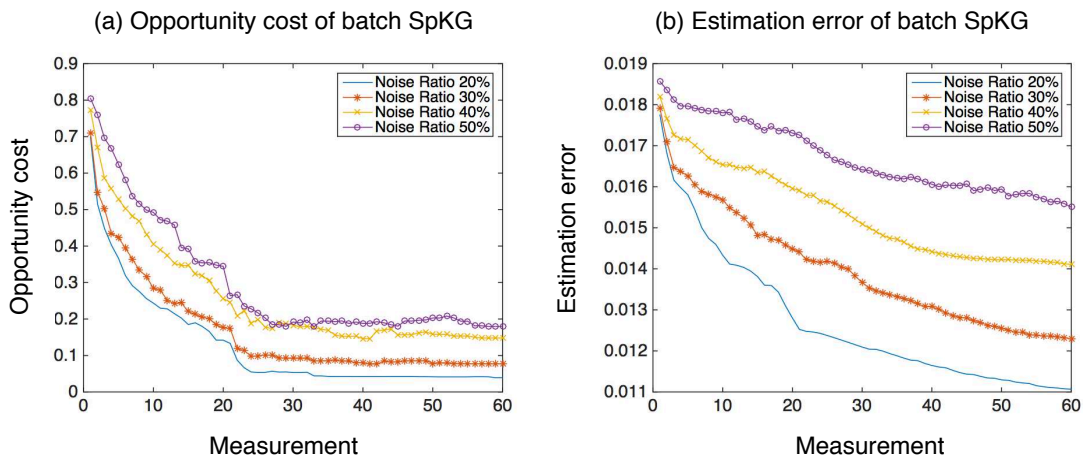}
\end{tabular}
\caption{(a) Averaged Percentage Opportunity Cost of Batch SpKG with Various Noise Ratios; (b) Averaged Normalized $\ell_2$ Estimation Error of Accessibility Coefficient with Various Noise Ratios}
\begin{flushleft}
\textit{Note.} Both figures are for replicated runs averaged over 200 times.
\end{flushleft}
\label{fig:4}
\end{figure*}

As one can see from Figure \ref{fig:3}, after five batch measurements, the estimate is still closer to the prior than the truth. After 10 batch measurements, we have discovered many of the accessible regions. After 15 batch measurements, we have not only discovered the location of the accessible regions, but also obtained good estimates for the actual accessibility value. And after 20 batch measurements, our estimate closely matches the truth. At last, we also try different noise levels $20\%, 30\%, 40\%$, and $50\%$ and repeatedly run such simulations for 200 times for each level. The averaged percentage OC and estimation error are plotted in Figure \ref{fig:4}. Here the normalized estimation error is the $\ell_2$ distance between the estimated coefficient and the truth divided by the length of the coefficient vector, which is 157 currently.  

Experimentally, fluorescent measurements are made by performing induction assays on prepared cell cultures. For each cell culture prepared, several samples are obtained, and fluorescence measurements are made via flow cytometry. Measurement dispersity in a small number of samples can be as large as $15\% \sim 20\%$ in standard deviation. For this noise level, we can see from the figure that most locations of the highly accessible regions can be found after about 25 observations. However, note that the true accessibility profiles sampled for the simulations are perturbed by a large amount. In reality, we suspect the truth to be more in agreement with the prior footprinting data, which implies better performance by the SpKG algorithm.

\subsubsection{Results using the Batch SpKG Algorithm with Length Mutagenesis}
In this set of experiments, we use the same set of probes in Section \ref{sec::simBKG} as the initial alternative library. The simulations are run with two different priors: a good prior and a bad prior. The bad prior is the one used in the above experiments, which is obtained by doing both vertical perturbation and horizontal shift from the in vitro DMS footprinting data. For the good prior, we only do vertical perturbation with the same amount, so we would think of the sparsity pattern of the good prior more proximal to the truth. Therefore, for this set of simulations with $L = 20$, $B = 3$, we set $w=10$ for the bad prior and $w=20$ for the good prior. 

\begin{figure*}[!htb]
\centering
\begin{tabular}{c}
\includegraphics[width=0.95\linewidth]{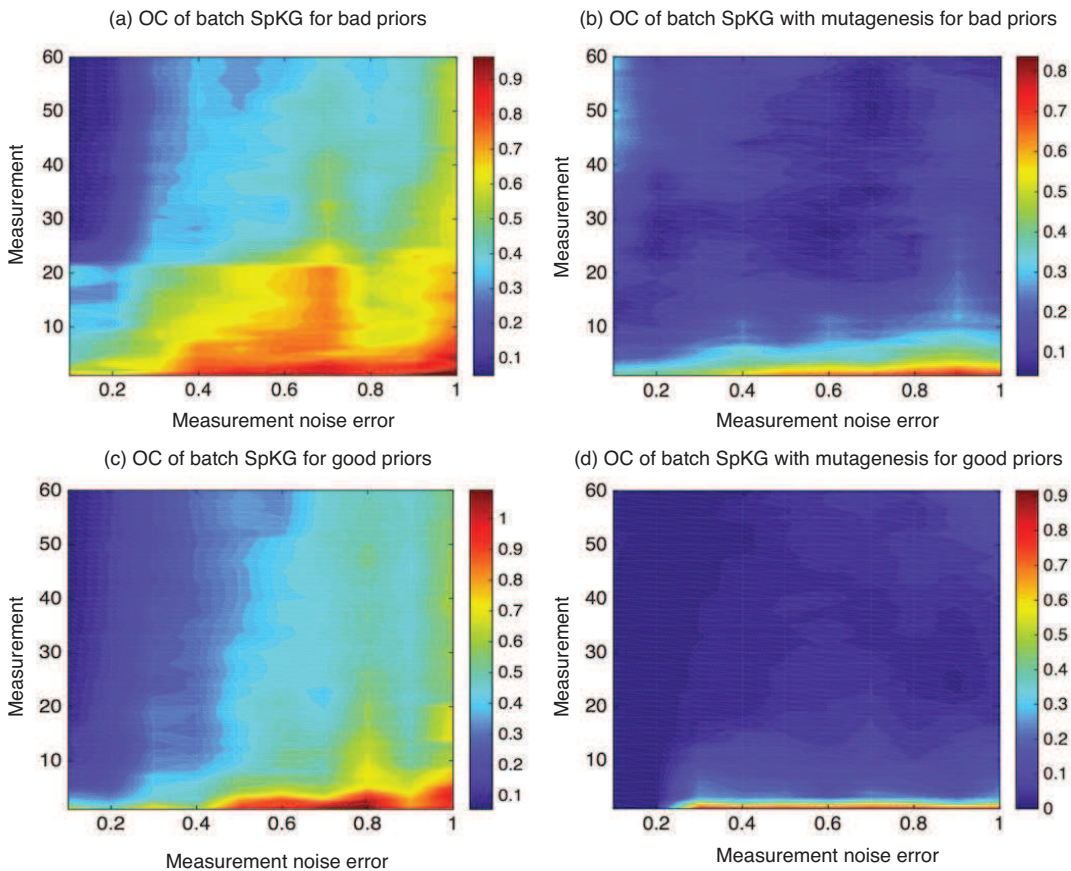}
\end{tabular}
\caption{Averaged Percentage Opportunity Cost of Batch SpKG with and without Mutagenesis for Good and Bad Priors over 300 Runs}
\begin{flushleft}
\textit{Note.} The contour plots show averaged percentage OC as a function of measurements and measurement noise errors. 
\end{flushleft}
\label{fig:6}
\end{figure*}

Figure \ref{fig:6} compares the batch SpKG algorithm with and without mutagenesis for both good and bad priors averaged over 300 simulation trials. It shows the mean percentage opportunity cost as a function of measurements and measurement noise errors. We try ten different noise levels from 10\% to 100\%.  For all of the figures, we observe that the OC\% decays as the number of measurements increases and the measurement noise decreases. Such plots can be useful in experimental budgeting and show the required number of measurements needed to obtain a certain level of optimality for some noise level. 

\begin{figure*}[!htb]
\centering
\begin{tabular}{c}
\includegraphics[width=0.95\linewidth]{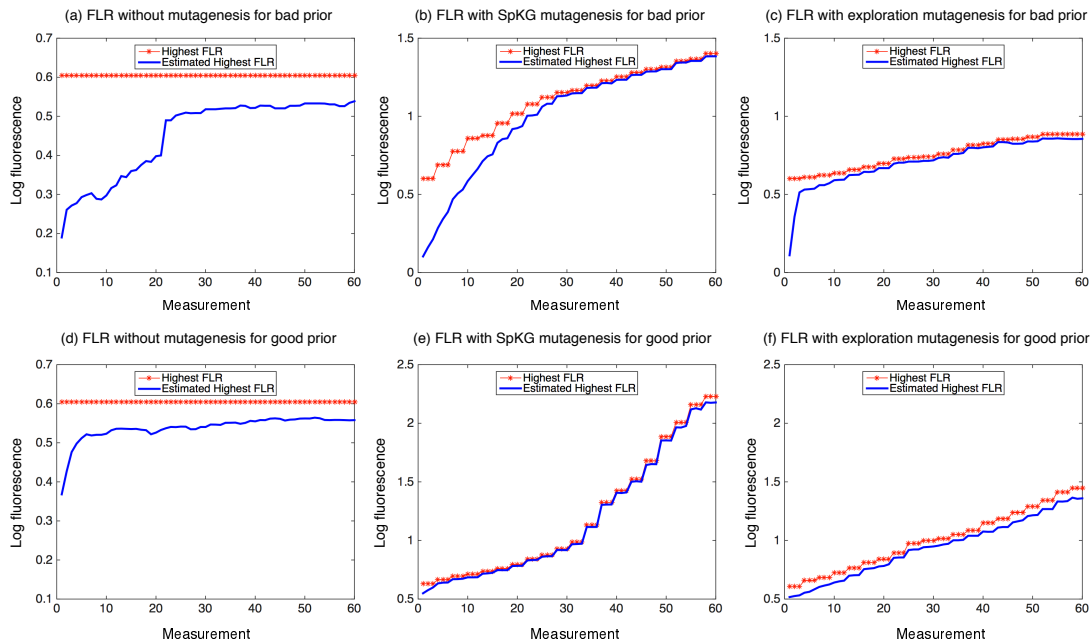}
\end{tabular}
\caption{True Highest Log Fluorescence and Estimated Highest Fluorescence for Good and Bad Priors Comparing Three Policies: Batch SpKG, Batch SpKG with Mutagenesis, and Exploration Mutagenesis}
\begin{flushleft}
\textit{Notes.} These plots show the values in true fluorescence the highest and the estimated highest. The difference between the two lines is the OC. The highest fluorescence for batch SpKG policy remains constant since we maintain the same probe alternative library throughout measurements. With mutagenesis, most of the time we could find the probe with higher fluorescence than before. For batch SpKG with mutagenesis, the new probe with the largest KG score is added. For exploration mutagenesis, the new probe is randomly added.  The results are averaged over 300 runs.
\end{flushleft}
\label{fig:8}
\end{figure*}

Comparing Figure \ref{fig:6}(a)(b) with (c)(d), we can see that the OC\% decreases much faster in the good prior cases as a function of the number of experiments, as expected. Comparing Figure \ref{fig:6} (a)(c) with (b)(d), we find that the OC\% generated from the  mutagenesis procedures tends to be smaller than those without mutagenesis. This is because when we add a new probe with the largest KG score into the alternative library with mutagenesis, we often add the one with a higher fluorescence level than the current fluorescences. In other words, the true highest fluorescence level is increasing as a new probe is added into the library. In such cases, it would be more intuitive to see how the highest fluorescence is varying over time. Figure \ref{fig:8} provides a more illustrative explanation for how the actual highest fluorescence changes. For this set of figures, we compare how the actual fluorescence changes over measurement. The red star line is the true highest fluorescence, and the blue solid line is the value of true fluorescence by estimation. So the difference between the two lines is the OC. We compare three different policies: batch SpKG, batch SpKG with mutagenesis, and exploration mutagenesis. Exploration mutagenesis involves randomly adding new probes if they are not in the current library. From Figure \ref{fig:8}, we can see that SpKG with mutagenesis has the ability to find new probes with fluorescence values about three to five times the highest values in the initial set. However, for exploration mutagenesis, the highest fluorescence improves less, as expected. This also proves the power of the KG policy to identify the potential alternative that can outperform the current optimal one. Furthermore, it is also worth noting that with mutagenesis, the value of the true fluorescence through the Bayesian estimation is pretty close to the truth. That means the OC is close to zero, which is consistent with Figure \ref{fig:6}(b)(d). 

\section{Conclusion}\label{sec::conclusion}

Identifying the accessibility pattern of an RNA molecule is an important topic in molecular biology. On one hand, the real experimental study is a long and expensive process for which adaptive learning procedures like the knowledge gradient policy are well suited to make experimental decisions. On the other hand, this problem naturally incorporates sparsity linear structure and thus requires more sophisticated statistical techniques to analyze the underlying model. To better help learn the accessibility profile of the RNA molecule, we use a recently derived SpKG policy \citep{li2015sparsekg}, which is a novel hybrid of Bayesian Ranking \& Selection and frequentist $\ell_1$ penalized regression approach called Lasso. This optimal learning algorithm has been shown to efficiently identify the accessibility pattern and learn the underlying sparsity structures. Controlled experiments also show that it outperforms several other policies.

Algorithmically, we also entend the SpKG policy into a general framework for batch mode learning, where the experimenter can run several experiments in batch. Empirical studies demonstrate the effectiveness of this policy for various experimental setups. Besides, we also derive the batch SpKG algorithm using length mutagenesis to expand the set of alternatives.
 In this procedure, the alternative library is adaptively enlarged as the most promising alternative is added in at each time. Controlled experiments also demonstrate its efficiency in identifying the accessibility pattern of the RNA molecule. 

In conclusion, it is worth noting that the SpKG algorithm has only been applied to sparse linear beliefs. Possible future directions of the work would include the study of more general nonlinear beliefs that incorporate sparsity structure. Despite this limitation,  we still believe the SpKG algorithm would allow efficient implementation for large data sets, and we would like to suggest this algorithm for solving more general application problems with sparse linear beliefs. 

\section*{Acknowledgement}
We thank Dr. Rick Russell for kindly providing the raw in vitro DMS footprinting data used in this work and previously published in \cite{russell2006paradoxical}.

% Appendix here
% Options are (1) APPENDIX (with or without general title) or 
%             (2) APPENDICES (if it has more than one unrelated sections)
% Outcomment the appropriate case if necessary
%
% \begin{APPENDIX}{<Title of the Appendix>}
% \end{APPENDIX}
%
%   or 
%
% \begin{APPENDICES}
% \section{<Title of Section A>}
% \section{<Title of Section B>}
% etc
% \end{APPENDICES}
%

% References here (outcomment the appropriate case) 

% CASE 1: BiBTeX used to constantly update the references 
%   (while the paper is being written).
\bibliographystyle{plainnat} % outcomment this and next line in Case 1
\bibliography{spkg} % if more than one, comma separated

% CASE 2: BiBTeX used to generate mypaper.bbl (to be further fine tuned)
%\input{mypaper.bbl} % outcomment this line in Case 2

\end{document}